\documentclass[11pt,a4paper]{article}

\usepackage{amssymb, amsmath, amsthm}

\newcommand{\E}{\mathbb{E}}
\newcommand{\p}{\mathbb{P}}
\newcommand{\R}{\mathbb{R}}

\newcommand{\bX}{\mathbf{X}}
\newcommand{\bA}{\mathbf{A}}
\newcommand{\bB}{\mathbf{B}}
\newcommand{\bnull}{\mathbf{0}}

\newcommand{\bk}{\mathbf{k}}

\newcommand{\bx}{\mathbf{x}}
\newcommand{\bY}{\mathbf{Y}}
\newcommand{\bZ}{\mathbf{Z}}

\newtheorem*{theorem}{Theorem}

\author{P\'eter Kevei\thanks{Bolyai Institute, University of Szeged, 
Aradi v\'ertan\'uk tere 1, 6720 Szeged, Hungary; 
e-mail: \texttt{kevei@math.u-szeged.hu}} \quad and \quad
P\'eter Wiandt\thanks{Bolyai Institute, University of Szeged, 
Aradi v\'ertan\'uk tere 1, 6720 Szeged, Hungary; 
e-mail: \texttt{wpeti88@gmail.com}}}

\title{Moments of the stationary distribution  of \\ subcritical 
multitype Galton--Watson processes \\ with immigration}

\date{}

\begin{document}

\maketitle

\begin{abstract}
In this short note we obtain necessary and sufficient
conditions for the existence of the moments of the 
stationary distribution of a subcritical multitype 
Galton--Watson process with immigration.

\end{abstract}

\section{Introduction}

Let $(\bX_n)_{n\geq 0} = (X_{n,1}, \ldots, X_{n,d})_{n \geq 0}$
be a $d$-type Galton--Watson process with immigration (GWI), 
defined as
\begin{equation} \label{eq:def-X}
\begin{split}
\bX_{n} & = \sum_{j=1}^d 
\sum_{i=1}^{X_{n-1,j}} \bA_{n,i;j} + \bB_{n}, \quad n = 1,\ldots, \\
\bX_0 & = \bnull = 
(0, \ldots,0 ),
\end{split}
\end{equation}
where $\bA_{n,i;j}, \bB_n$, $n \geq 1$, $i \geq 1$, $j \in \{1,\ldots ,d\}$,
are independent $d$-dimensional
random vectors with integer coordinates such that
$\{ \bA_{n,i;j}: n \geq 1, i \geq 1 \}$ is an identically distributed 
sequence of random variables for each $j \in \{1, \ldots, d\}$, 
and $\bB_n$,  $n\geq 1$, are identically 
distributed. Here $X_{n,j}$ is the number of $j$-type individuals 
in generation $n$, $\bA_{n,i;j}$ are the offsprings produced by
the $i$th individual of type $j$ in generation $n-1$, and 
$\bB_n$ are the immigrants. In what follows, vectors,
both deterministic and random, are denoted 
by boldface letters, and are meant as $d$-dimensional row vectors.

\smallskip

Branching processes play an important role in models of 
genetics, molecular biology, physics and computer science. 
As a main reference on branching processes
we refer to the classical books by Athreya and Ney \cite{AthreyaNey},
by Mode \cite{mode} and by Haccou et al.~\cite{Haccou}.

Multitype Galton--Watson processes with immigration was
introduced and studied by Quine \cite{Quine}. 
In \cite{Quine} necessary and sufficient condition was 
given for the existence of stationary distribution in the subcritical case.
A complete answer for the existence of a limiting stationary distribution
for multitype GWI was obtained by Kaplan \cite{kaplan}. Also Mode \cite{mode} 
gave a sufficient condition for the existence of a stationary distribution.
Multitype GWI processes in random environment were investigated
by Roitershtein \cite{Roit07}, Roitershtein and Zhong \cite{RoitZhong},
Wang \cite{Wang}, to mention just a few.

However, interestingly enough, conditions for the existence of moments
of the stationary distribution are not known except in special cases.
Explicit formula for the variance was obtained by 
Quine \cite{Quine} and for the third moment very recently by Barczy et 
al.~\cite[Lemma 1]{Barczy}.
In the unpublished work by Sz\H{u}cs \cite{Szucs} under 
ergodicity assumptions the existence of general moments of the 
stationary distribution were obtained for multitype GWI processes.
Even in the single type case, we are only aware of the recent results
by Buraczewski and Dyszewski \cite{BD} on Galton--Watson processes
in random environment (without immigration), from which 
the existence of certain moments can be deduced, 
see Lemma 3.1 in \cite{BD}, and by 
Basrak and Kevei \cite[Lemma 1]{Basrak}
on GWI processes in random environment.

The aim of the present paper is to fill this gap and 
obtain necessary and sufficient conditions for the 
existence of the moments of the stationary distribution.

\section{Results}

We assume that the offspring means are finite,
and let $M$ denote the offspring mean matrix
\begin{equation} \label{eq:defM}
M = 
\begin{pmatrix}
\E \bA_{1, 1; 1} \\
\vdots \\
\E \bA_{1, 1; d} \\
\end{pmatrix}
= 
\begin{pmatrix}
m_{1,1} & \ldots & m_{1,d} \\
\vdots & \ddots & \vdots \\
m_{d, 1} & \ldots & m_{d,d} \\
\end{pmatrix},
\end{equation}
that is $m_{i,j}$ is the mean offsprings of type $j$ produced by an 
individual of type $i$. Let $\rho$ denote the spectral radius of
$M$, and assume that the process is subcritical, i.e.~$\rho < 1$.
Without immigration a subcritical process dies out almost surely
exponentially fast. While, Quine \cite{Quine} showed that 
immigration with finite logarithmic moments ensures the existence 
of a stationary distribution.
To ease notation we introduce the random operators $\theta_n$ as
\begin{equation} \label{eq:def-theta}
\theta_n \circ \bk =  
\sum_{j=1}^d \sum_{i=1}^{k_j} \bA_{n, i; j},
\quad \bk = (k_1, \ldots, k_d).
\end{equation}
We slightly abuse the notation writing
$\theta_n \circ (\bk_1 + \bk_2) = \theta_n \circ \bk_1 + \theta_n \circ \bk_2$,
where on the right-hand side the two summands are independent.
Further, write $\Pi_n = \theta_1 \circ \ldots \circ \theta_n$, for $n \geq 1$,
and $\Pi_0 = \mathrm{Id}$.
With this notation \eqref{eq:def-X} can be written as
$\bX_n = \theta_n \circ \bX_{n-1} + \bB_n$. Iteration gives 
that the stationary distribution can be represented 
in distribution as
\begin{equation} \label{eq:Y-stat}
\bY =
\bB_1+ \theta_1 \circ \bB_2+ \theta_1 \circ \theta_2 \circ \bB_3 + \ldots
=\sum_{i=0}^{\infty} \Pi_i \circ \bB_{i+1}, 
\end{equation}
provided that the infinite sum exists.
This corresponds to formula (16) in \cite{Quine} in terms of generating
functions.

Introduce the notation
$\| \bx \| = \sum_{j=1}^d |x_j|$ for the $\ell^1$ norm of a vector
in $\R^d$, and also for the generated matrix norm 
$\| A \| = \sup_{\| \bx \| = 1} |\bx A|$. Note that $\| A \|$ is the maximum
absolute row sum, since we multiply from the left.

Our main result is the following.
\begin{theorem} \label{thm:main}
Let $(\bX_n)_{n \geq 0}$ be a $d$-type, subcritical 
Galton--Watson process with immigration.
If $\E \| \bA_{1,1;i} \|^{\max\{\alpha, 1\}} < \infty$ for all 
$i \in \{1, \ldots, d\}$, and $\E \| \bB \|^\alpha < \infty$ for some 
$\alpha > 0$, then 
$\E \| \bY \|^\alpha < \infty$. In particular, each component of $\bY$ has 
finite moment of order $\alpha$.
\end{theorem}

Note that by the subcriticality assumption $\E \| \bA_{1,1;i} \| < \infty$
for any $i \in \{1,\ldots, d\}$. Thus, even for $\alpha < 1$ we assume 
the existence of the offspring mean. However, the immigration distribution 
might have infinite mean.

Also note that we do not assume the existence of a unique stationary 
distribution. We only use that the conditions of the theorem implies 
that the infinite sum in \eqref{eq:Y-stat} is a.s.~finite.

Finally, we mention that these are only sufficient conditions.
Clearly, if a type-1 particle never immigrates, then the 
distribution of $\bA_{1,1;1}$ does not matter in the 
stationary distribution. Otherwise, from formula \eqref{eq:Y-stat}
we see that our conditions are also necessary.

\section{Proof}

\subsection{The case $\alpha \geq 1$}

Let $\alpha \geq 1$ be fixed.
First we prove the theorem under the additional assumption that
\begin{equation} \label{eq:ass-1}
\| M \| < 1,
\end{equation}
that is the row sums are less than 1.
Let $\mu_j = \sum_{\ell=1}^d m_{j,\ell}$, $j \in \{1, \ldots, d\}$.

We prove that  
\begin{equation} \label{eq:def-Mk}
M_\alpha (k)=
\E \|\Pi_k \circ \bB_{k+1} \|^\alpha
\end{equation}
decreases exponentially.

Recall \eqref{eq:def-theta}. To ease notation we suppress the lower index.
We have 
\[ 
\| \theta \circ \bk \| = \sum_{j=1}^d \sum_{i=1}^{k_j} \| \bA_{i;j} \|
= \sum_{j=1}^d \sum_{i=1}^{k_j} \sum_{\ell=1}^d A_{i;j}^{(\ell)}
=: \sum_{j=1}^d S_{k_j;j},
\] 
where $\| \bA_{i;j} \| = \sum_{\ell=1}^d A_{i;j}^{(\ell)}$ is 
the number of all the offsprings of the $i$th individual of type $j$,
with $A_{i;j}^{(\ell)}$ being the number of $\ell$-type offsprings,
and $S_{k;j} = \sum_{i=1}^k \| \bA_{i;j} \|$ is the sum of $k$ iid 
\emph{scalar} random variables, with mean
$\E \| \bA_{1;j} \| = \sum_{\ell=1}^d m_{j,\ell} = \mu_j < 1$.

For any  $\overline \mu \in (\max_{1 \leq j \leq d} \mu_j, 1)$
there exists $k_0'$ such that 
\begin{equation} \label{eq:k0}
\E \left( \frac{S_{k;j}}{k} \right)^\alpha  < \overline \mu
\quad \text{for all } k \geq k_0', \ j \in \{1, \ldots, d \}.
\end{equation}
This follows from the strong law of large numbers combined 
with the uniform integrability of $S_{k;j}/k$. Put
\begin{equation} \label{eq:c0}
\max_{1 \leq j \leq d } \max_{ k < k_0'} 
\E \left( \frac{S_{k;j}}{k} \right)^\alpha = c_0.
\end{equation}
Since the function $x^\alpha$ is convex, we have
\[
\begin{split}
\left( 
\frac{\sum_{j=1}^{d} S_{k_j; j}}{k_1 + \ldots + k_d} 
\right)^\alpha  
& = \left( \sum_{j=1}^{d}  \frac{k_j}{k_1 + \ldots + k_d} 
\frac{S_{k_j; j}}{k_j} 
\right)^\alpha  \\
& \leq \sum_{j=1}^{d}  \frac{k_j}{k_1 + \ldots + k_d} 
\left( \frac{S_{k_j; j}}{k_j} \right)^\alpha.
\end{split}
\]
Combined with \eqref{eq:k0} and \eqref{eq:c0} this implies
\[
\begin{split}
\E \| \theta \circ \bk \|^\alpha & 
=  \| \bk \|^\alpha \E 
\left( \frac{\| \theta \circ \bk \|}{\|\bk \|}\right)^\alpha \\
& \leq \| \bk \|^\alpha \sum_{j=1}^d \frac{k_j}{\| \bk \| }
\E \left( \frac{S_{k_j;j}}{k_j} \right)^\alpha \\
& \leq \| \bk \|^\alpha \sum_{j=1}^d 
\frac{k_j}{\| \bk \| } 
\left[ I( k_j \geq k_0') \overline \mu + I(k_j < k_0') c_0 \right] \\
& \leq \| \bk \|^\alpha \left( \overline \mu + \frac{k_0' c_0 d}{\| \bk \|}
\right),
\end{split}
\]
where $I(\cdot)$ stands for the indicator function.
Choosing $k_0 > 2 k_0' c_0 d / ( 1 - \overline \mu)$ we obtain that
\begin{equation} \label{eq:aux1}
\E \| \theta \circ \bk \|^\alpha \leq \mu \| \bk \|^\alpha
\quad \text{whenever} \ \| \bk \| \geq k_0, 
\end{equation}
with $\mu = ( 1 + \overline \mu)/2 < 1$.

Let $\bZ$ be a random vector with nonnegative integer components,
independent of the $\bA$'s. 
Put 
\[
c_1 = \max_{|\bk| < k_0} \E  \| \theta \circ \bk \|^\alpha.
\] 
Then, by \eqref{eq:aux1}
\[
\begin{split}
\E \| \theta \circ \bZ \|^\alpha 
& = \sum_{\bk} \p ( \bZ = \bk ) \E \| \theta \circ \bk \|^\alpha \\
& = \sum_{\bk: \| \bk \| \geq k_0} 
\p ( \bZ = \bk ) \E \| \theta \circ \bk \|^\alpha + 
\sum_{\bk: 0 < \| \bk \| < k_0} 
\p ( \bZ = \bk ) \E \| \theta \circ \bk \|^\alpha \\
& \leq \sum_{\bk: \| \bk \| \geq k_0} 
\p ( \bZ = \bk ) \mu \| \bk \|^\alpha + 
\sum_{\bk: 0 < \| \bk \| < k_0} 
\p ( \bZ = \bk ) c_1 \\
& \leq \mu \E \| \bZ\|^\alpha + c_1 \, \E \| \bZ \|.
\end{split}
\]
Turning back to $M_\alpha(k)$ in \eqref{eq:def-Mk}, we obtain the 
recursion
\begin{equation} \label{eq:M-recursion}
M_\alpha(k) \leq \mu M_\alpha(k-1) + c_1 M_1(k-1), 
\end{equation}
with $\mu < 1$. Note that $M_1(k)$ is the expectation of 
the total number of individuals in generation $k$ in a multitype 
Galton--Watson process without immigration, starting with 
$\bB_{k+1}$ at generation 0. Therefore 
$M_1(k) = \| M^k \E \bB \|$, which decreases exponentially fast.
Thus, recursion \eqref{eq:M-recursion} implies that 
for some $\nu \in (\mu, 1)$ and $C > 0$
\[ 
M_\alpha (k) \leq C \nu^k \quad \text{for all } k \geq 0.
\] 

The statement now follows from Minkowski's inequality, as
\[
\begin{split}
\left( \E \| \bY \|^\alpha \right)^{1/\alpha} & = 
\left( \E \Big\| \sum_{k=0}^\infty \Pi_k \circ B_{k+1} \Big\|^\alpha
\right)^{1/\alpha} \\
& \leq \left( \sum_{k=0}^\infty M_\alpha (k) \right)^{1/\alpha} < \infty.
\end{split}
\]

\bigskip

The additional assumption in \eqref{eq:ass-1} can be omitted easily.
By Gelfand's formula for the spectral radius we have
\[
\lim_{k \to \infty} \| M^k \|^{1/k} = \rho,
\]
which is strictly less than 1, by subcriticality. Thus, there 
exists $r$ such that $\| M^r \| < 1$. The matrix $M^r$ is the mean 
matrix of the offspring distribution corresponding to $\Pi_r$, 
i.e.~when we sample the process only in every $r$th generation.
Therefore, the previous argument gives that 
$M_\alpha(rk + i)$ in \eqref{eq:defM} decreases exponentially
for each $i \in \{0,1, \ldots, r-1\}$.
Clearly, then $M_\alpha(k)$ also decreases exponentially, and the 
result follows from Minkowski's inequality as above.

\subsection{The case $\alpha < 1$}

This case is in fact simpler, but needs to be treated differently.

First, assume again that \eqref{eq:ass-1} holds.
We use the same notations as above. Now $x^\alpha$ is concave, thus
by Jensen's inequality for any $\bk$
\[
\begin{split}
\E \| \theta \circ \bk \|^\alpha 
& = \E \big( \sum_{j=1}^d S_{k_j;j} \big)^\alpha \leq 
\Big( \E \sum_{j=1}^d S_{k_j;j} \Big)^\alpha \\
& = \left( \sum_{j=1}^d k_j \mu_j \right)^\alpha
\leq \mu \| \bk \|^\alpha,
\end{split}
\]
with $\mu = \max_{1 \leq j \leq d } \mu_j$, implying
\[
\begin{split}
\E \| \theta \circ \bZ \| 
& = \sum_{\bk} \p ( \bZ = \bk ) \E \| \theta \circ \bk \|^\alpha \\
& \leq \sum_{\bk} \p ( \bZ = \bk ) \mu \| \bk \|^\alpha  
\leq \mu \E \| \bZ\|^\alpha.
\end{split}
\]
Therefore, the exponential decrease of $M_\alpha(k)$ in \eqref{eq:def-Mk}
follows. By subadditivity we have
\[
\E \| \bY \|^\alpha  = 
\E \Big\| \sum_{k=0}^\infty \Pi_k \circ B_{k+1} \Big\|^\alpha
\leq \sum_{k=0}^\infty M_\alpha (k)  < \infty,
\]
as claimed.

Condition \eqref{eq:ass-1} can be omitted the same way as before.

\bigskip
\noindent \textbf{Acknowledgement.}
P\'eter Kevei is supported by the J\'anos Bolyai Research Scholarship of the Hungarian
Academy of Sciences,  by the NKFIH grant FK124141, and 
by the EU-funded Hungarian grant EFOP-3.6.1-16-2016-00008.
P\'eter Wiandt is supported  
by the Ministry of Innovation and Technology grant UNKP-19-2-SZTE-139,
by the Ministry of Human Capacities, Hungary grant 20391-3/2018/FEKUSTRAT,
and by the EU-funded Hungarian grant EFOP-3.6.2-16-2017-00015.


\begin{thebibliography}{10}

\bibitem{AthreyaNey}
K.~B. Athreya and P.~E. Ney.
\newblock {\em Branching processes}.
\newblock Springer-Verlag, New York-Heidelberg, 1972.
\newblock Die Grundlehren der mathematischen Wissenschaften, Band 196.

\bibitem{Barczy}
M.~Barczy, F.~K. Ned\'{e}nyi, and G.~Pap.
\newblock On aggregation of multitype {G}alton-{W}atson branching processes
  with immigration.
\newblock {\em Mod. Stoch. Theory Appl.}, 5(1):53--79, 2018.

\bibitem{Basrak}
B.~Basrak and P.~Kevei.
\newblock Limit theorems for branching processes with immigration in a random
  environment.
\newblock https://arxiv.org/abs/2002.0063, 2020.

\bibitem{BD}
D.~Buraczewski and P.~Dyszewski.
\newblock Precise large deviation estimates for branching process in random
  environment.
\newblock https://arxiv.org/abs/1706.03874v1, 2019.

\bibitem{Haccou}
P.~Haccou, P.~Jagers, and V.~A. Vatutin.
\newblock {\em Branching processes: variation, growth, and extinction of
  populations}, volume~5 of {\em Cambridge Studies in Adaptive Dynamics}.
\newblock Cambridge University Press, Cambridge; IIASA, Laxenburg, 2007.

\bibitem{kaplan}
N.~Kaplan.
\newblock The multitype {G}alton-{W}atson process with immigration.
\newblock {\em Ann. Probability}, 1(6):947--953, 1973.

\bibitem{mode}
C.~J. Mode.
\newblock {\em Multitype branching processes. {T}heory and applications}.
\newblock Modern Analytic and Computational Methods in Science and Mathematics,
  No. 34. American Elsevier Publishing Co., Inc., New York, 1971.

\bibitem{Quine}
M.~P. Quine.
\newblock The multi-type {G}alton-{W}atson process with immigration.
\newblock {\em J. Appl. Probability}, 7:411--422, 1970.

\bibitem{Roit07}
A.~Roitershtein.
\newblock A note on multitype branching processes with immigration in a random
  environment.
\newblock {\em Ann. Probab.}, 35(4):1573--1592, 2007.

\bibitem{RoitZhong}
A.~Roitershtein and Z.~Zhong.
\newblock On random coefficient {INAR}(1) processes.
\newblock {\em Sci. China Math.}, 56(1):177--200, 2013.

\bibitem{Szucs}
G.~Sz\H{u}cs.
\newblock Ergodic properties of subcritical multitype {G}alton-{W}atson
  processes.
\newblock Available on arXiv: https://arxiv.org/abs/1402.5539.

\bibitem{Wang}
H.~M. Wang.
\newblock A note on multitype branching process with bounded immigration in
  random environment.
\newblock {\em Acta Math. Sin. (Engl. Ser.)}, 29(6):1095--1110, 2013.

\end{thebibliography}
\end{document}